\documentclass[a4paper, 10pt, conference, onecolumn]{ieeeconf}      

\IEEEoverridecommandlockouts                              
\def\ignore#1{}
\usepackage{epsfig} 
\usepackage{psfrag}
\usepackage{amsmath} 
\usepackage{amssymb}  

\title{\LARGE \bf Optimal control theory~:\\ a method 
for the design of wind instruments}

\author{G. Le Vey 
\thanks{IRCCyN UMR-CNRS 6597 and Ecole des Mines de Nantes - 4, rue A. Kastler - BP 20722 - 44307 NANTES Cedex 3, France.
{\tt\small levey@emn.fr}}}
\begin{document}
\maketitle
\thispagestyle{empty}
\pagestyle{empty}
\bibliographystyle{plain}
%
%
\begin{abstract}
It has been asserted previously by the author that optimal control 
theory can be a valuable framework for theoretical studies about the shape that a wind instrument should have 
in order to satisfy some optimization criterion, inside a fairly 
general class. The purpose of the present work is to develop 
this new approach with a look at a specific criterion to be optimized. 
In this setting, the Webster horn equation is regarded as a controlled dynamical equation in the space variable. 
Pressure is the state, the control being made of two parts~: one variable part, the inside diameter of the 
duct and one constant part, the weights of the elementary time-harmonic components of the 
velocity potential. Then 
one looks for a control that optimizes a criterion related to the definition of an {oscillation regime} 
as the cooperation of several natural modes of 
vibration with the excitation, the {playing frequency} being the one that maximizes the total generation of energy, as exposed by A.H. Benade, following H. Bouasse. At the same time the relevance of this criterion is questionned with 
the simulation results.
\end{abstract}
%
\section{Introduction}
\label{intro}
Designing high quality musical instruments has been the main concern of makers for centuries, accumulating know-how from their 
predecessors by some kind of trial-and-error process, involving musicians. On another hand, since Bernouilli and Lagrange 
in the XVIII$^{th}$ century and the important writings of Helmholtz in the XIX$^{th}$ century a great deal of research has also been conducted in order to understand 
the principles underlying sound production in traditional or modern musical instruments~: see e.g.~\cite{benade:1990,bouasse:1929,campbell:1987,kergomard:2008,fletcher:1991} for presentations of this vast subject.  
Wind instruments such as brasses and woodwinds, making use of an air column, must be designed in such a way as to properly arrange the natural frequencies of this air column, 
in order that {\em regimes of oscillation} can take place in conjunction with the nonlinear flow-control device (reed or player's lips)~\cite{benade:1990}. 
As a consequence, harmonicity requirements for wind instruments between natural frequencies are desirable properties, due to the fact that 
intonation, responsiveness and tonal colour can be attributed to well established physical 
properties of the instrument~\cite{summermeeting:2009}. 
It is generally accepted that these harmonicity requirements can lead to rather different shapes, 
which can be (piecewise) cones, cylinders or more complex, as can can be seen from 
the actual duct of real wind instruments. 
Nevertheless, see~\cite{dalmont-kergomard:1994} for interesting results about piecewise conic ducts, 
following an approach based on characteristic impedances and 
transfer matrices. Thus there is room for developping new methods that 
could help for the design. The present work is dedicated to such a task, with main accent 
on the methodological aspects.\\ 
In~\cite{kausel:2001}, transmission line modelling of horns is used 
together with finite dimensional 
optimization techniques and a practical tool~\cite{kausel:1999} 
was developed for instrument makers, 
in order to help improving existing instruments as
well as to design new ones according to a given specification.
Desired properties, specified in musical 
terms like intonation, response and pitch variability were compared with 
calculated values based on an instrument's actual geometry. 
Similarly, in~\cite{noreland:2003-2}, a method is presented for 
optimizing the shape of a brass instrument with respect to its 
intonation and impedance peak magnitudes. 
The instrument is modelled using a one-dimensional transmission 
line analogy with truncated cones. Through the use of an appropriate 
choice of design variables, 
the finite dimensional optimization finds smooth horn profiles, 
that can also help 
in correcting existing instruments, the shape having been 
designed from an a priori choice 
of a succession of cones. 
In~\cite{forbes:2006}, a frequency-domain method, using inverse quantum scattering for the one-dimensional 
Klein-Gordon equation, allows to recover the area function of a given acoustical duct in a noninvasive way, 
without measuring directly neither the input impedance nor the reflectance~: 
this last one is mathematically derived from the wave radiated in response to a high-impedance source.
In~\cite{inacio:2007}, parameter optimization techniques have been used to design shapes for brass-trapping 
Helmholtz resonators that resonate at a design set of acoustic eigenvalues, taking into account 
physical and geometrical constraints. 
In~\cite{henrique:2003}, a finite-element eigenanalysis model of bars is used to compute optimal shapes 
for mallet percussion instruments such as vibraphone or marimba-type bars. 
The objective function for the optimization procedure is a targeted set of modal frequencies.\\
As one can see, optimization is at the heart of research on the design question for musical instruments and 
the eigenstructure often plays a central role in the optimization criteria. 
The shape optimization of wind instruments, an {inverse problem}, appears to be a complex one, 
when compared to the analysis of what physically occurs inside a given 
instrument, a direct, somewhat easier albeit not at all obvious, 
{problem}, as can be seen from the great amount of research that it is the subject of. 
The purpose of the present work is a continuation of preliminary ones by the 
author~\cite{levey-isma07,levey-med08}, where the design question for axisymmetric wind instruments was posed as 
an optimal control problem for the so-called {\em Webster equation} through a suitable 
reformulation as a dynamical equation in the spatial dimension, 
thereby using infinite dimensional optimization of optimal control 
theory for continuous systems. This is in sharp contrast with the above mentionned references which 
use finite-dimensional optimization schemes, i.e. with a finite number of parameters as unknowns in the theoretical 
formulation. 
For more details about the horn equation, see~\cite{eisner:1967} or the more recent~\cite{rienstra:2005} 
for a more mathematical and updated exposition. 
Here, besides the exposition of how optimal control theory fits the wind instruments design question, 
a special attention is given to the choice of one specific 
criterion, on the basis of physical considerations, whereas in~\cite{levey-isma07}, only a fairly general class of criteria -affine in the control- were considered that led 
to a general result about the structure of the shape for the duct of a wind instrument. \\
The chosen criterion is in no way unique~: many other choices could be made, 
depending on which characteristics (physical, perceptive and so on) are focused on. 
It appeared to be interesting to investigate from a musical acoustics point of view as it comes from 
one definition of an {\em oscillation regime} inside wind instruments, understood as the cooperation 
of several natural modes of vibration, 
the 'playing frequency' being the one that maximizes the total generation of 
energy~\cite{bouasse:1929,benade:1990}. At this stage, it has some 
arbitrariness in it~: the physical relevance of this criterion is not taken here for granted, as 
no other published material has been found (than that of A.H. Benade) using this notion of oscillation regime.
The reason for the choice of this criterion 
is twofold : first, to illustrate the general methodology that uses optimal control 
theory ; second, to question the relevance of this criterion, that appeared to be appealing 
at first sight. 
The purpose then is to question how interesting it is for designing wind instruments more than 
to assert such a relevance a priori~: 
numerical results in section~\ref{numeric} show some qualitative features linked to the chosen criterion. \\
Nevertheless, on the basis of this criterion, in the context of optimal control theory, 
necessary conditions coming from the strong Pontryagin Maximum Principle~\cite{pontriaguine:1978,clarke:1983} 
are derived and some numerical simulations are conducted, 
in order to illustrate the theoretical predictions. 
In the presented approach, the Webster horn equation is regarded as a controlled {dynamical system} in one spatial variable (the axis of an axisymmetric instrument). 
The pressure (or the {velocity potential}) is regarded as the {state} and the {control} is made of two parts~: one is the variable diameter 
(actually its derivative with respect to the axial space dimension) of the duct. The second part of the control is 
made of the constant weights of a time-harmonic decomposition of the velocity potential, which are 
unknown, according to the chosen criterion, the energy inside the duct. 
Then one looks for a control that maximizes the latter. 
This is another difference with the approach in~\cite{kausel:2001,noreland:2003-2}, where the measure is taken to be the deviation from a measured input impedance peaks location, 
making use of a {model-fitting oriented} approach. Whereas these works can be described as inverse problems, the present 
one falls into the category of design problems. But this distinction is a matter of convenience as both problems 
only differ from each other by the chosen criterion, the unknown being the duct shape in each case. 
It is worth noting in that respect that a similar criterion 
as in~\cite{noreland:2003-2} could be chosen here as well, while defining the corresponding output of the 
model, leading to the same problem of fitting a model to given data. This will be the subject 
of future research, with interesting comparisons to be made, the main point here being methodological. 
The approach followed hereafter shows that, with the chosen optimization 
criterion and the modelling approach, the obtained shape is piecewise continuous. 
Moreover, when taking into account more and 
more time-harmonic components, is appears from the simulations that the the obtained duct shape is a cone~: 
such a result tends to confirm that this criterion is likely not to be sufficient to design high-quality instruments 
but the realistic mode shapes that are obtained are encouraging. 
It should be noted that, to the best knowledge of the author, 
this way of looking at the classical {Webster horn equation}, as a controlled dynamical equation, appears to be new. 
One important point is that extensions 
of the optimal control methodology to other more realistic, nonlinear models of 
wind instruments and design criteria can be considered, although Fourier analysis will generally be no more 
valid in such contexts~: the choice of an adequate criterion will remain crucial. 
Numerical simulations show that the proposed approach is appealing for design purposes. 
Nevertheless, the presented work has to be considered as a preliminary study that opens new perspectives for the 
design of wind instruments~: it is not to be understood as 
an achievement in itself and has surely to be investigated more deeply for real applications. In particular, 
several aspects are not touched here such as taking account of
tone holes in woodwinds (flutes, oboes, etc..), incorporating the nonlinear excitation mechanism 
(reed, lips for brasses) or 
choosing other criteria, in order to include perceptive parameters (such as given by a 
musician), etc...Much remains to be done.\\
The paper is organized as follows~: section~\ref{preliminaries}, in order to fix the notations, 
recalls the standard horn equation and the 
definition of an oscillation regime as used here. 
Basic results from Optimal Control theory are briefly given and lead to the formulation of the horn 
equation as a controlled dynamical system. 
Then a detailed exposition of how the oscillation regime is formulated is given in section~\ref{criterion}. 
Section~\ref{webster} gives the mathematical formulation of the optimization problem to solve, through 
the necessary conditions 
coming from optimal control theory. 
Numerical results are given in section~\ref{numeric}~: they have to be considered as preliminary results 
as no thorough or exhaustive investigation have been 
conducted at this stage. Eventually, some conclusions and perspectives are drawn in section~\ref{conclusion}.
An appendix provides a variational derivation of the horn wave equation and 
reminds basic results from Sturm-Liouville theory. 
\section {Preliminaries~: Horn wave equation, Oscillation regime, Optimal control theory}
\label{preliminaries}
\begin{figure}[h]
\psfrag{D}{$D$}
\psfrag{DD}{$D(x)$}
\psfrag{xx}{$x$}
\psfrag{OO}{$0$}
\psfrag{L}{$L$}
\centerline{\psfig{figure=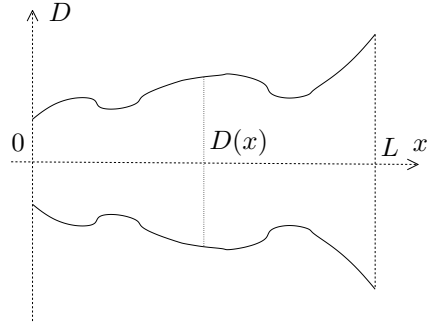,width=5.5cm}}
\caption{Geometry of an axisymmetric horn}\label{fig1}
\end{figure}
\subsection{\underline{Wave horn equation, oscillation regime}}
\label{hornoscillation}
Consider an axisymmetric horn with length $L$ and section diameter $D(x)$, a function of the 
independent variable, the space dimension $x$, according to Fig. 1. Let 
$\rho_0$ be the medium mass density, $p(x,t) $ the acoustic pressure in the 
medium, i.e. the deviation from the atmospheric pressure, 
and $v(x,t)$ the {particular velocity} inside the horn at abcissa $x$.
Let $k=\omega/c$ be the wave number, $c$ being the sound velocity 
and $\omega$ the angular frequency 
in case of a harmonic regime that is focused on here. 
For derivatives 
of any quantity the index notation will be used in this section~: for example the time (resp. 
space) first derivative of $\phi$ will be denoted~ $\phi_t$ (resp. $\phi_x$). 
Let us first recall some usual modelling hypothesis~: i) the diameter $D(x)$ is a {slowly} varying function of $x$~; ii) acoustical quantities are 
functions of $x$ and $t$ only~; iii) a {velocity potential} $\phi$ exists, which means that $v=\nabla\phi=\phi_x$~; iv) one focuses on {standing wave, time-harmonic solutions} that 
write~: $\phi(x,t)=\Re(\varphi(x)e^{i\omega t})$ but one will omit the $\Re$ part in the sequel as it is related 
to the time dependence that will not be involved in the optimization process below~; 
v) in this work there is no dissipation. 
It can be shown that, in this situation, the pressure at $x$ writes~: $p= - \rho_0 \phi_t$. 
It is worth noticing that the first hypothesis is necessary for the 1D plane wave, 
Webster equation to be valid~\cite{inacio:2007} or 
\cite{kergomard:2008}, chap.~6. 
Another limitation to this equation is to lower frequency modes. In~\cite{kergomard:2008}, the following 
approximate condition on the radius derivative $R^{'}$ is given for the horn equation to be valid for plane waves~:
\begin{equation}
\frac{1}{2}\int_0^Lk{R^{'}}^2(x) dx \ll 1
\label{eq.inequal}
\end{equation}
This condition will be used to fix bounds on the control in the simulations presented in section~\ref{numeric}. 
With these conventions, the {Webster horn equation} is easily obtained as 
the Euler equation of a variational problem~:
\begin{eqnarray}
\frac{1}{c^2}\phi_{tt}-\phi_{xx}-2\frac{D_x}{D}\phi_x=0
\label{eq.euler}
\end{eqnarray}
which reduces to~:
\begin{equation}
\begin{array}{l}
\phi^{''}+2\frac{D^{'}}{D}\phi^{'}+k^2 \phi=0
\label{eq.sturmliouville}
\end{array}
\end{equation}
in case of a harmonic regime of interest here (see appendix a for details). Boundary conditions are then 
given for this equation, depending on the physical situation under study, making it a Sturm-Liouville 
problem (see appendix b for a brief reminder of two useful results).\\
On another side, the notion of {\em oscillation regime} can be defined, 
following H. Bouasse~\cite{bouasse:1929} and A.H. Benade~\cite{benade:1990}. 
It is {\em ``that state of the 
collective motion of an air column in which 
a nonlinear excitation mechanism collaborates with a set of air column modes to 
maintain a steady oscillation containing several harmonically related frequency components, 
each with its own definite amplitude. Then, the 'playing frequency' is the one that 
maximizes the total generation of energy''}~\cite{benade:1990}. \\
Chosing such a definition as the basis of an optimization criterion 
relies upon physical/acoustical considerations of what happens inside the instrument. Nevertheless, 
it has some arbitrariness when considering perceptual quality of sound, as discussed in the introduction. 
But further investigations looked desirable. Other criteria will be considered 
in future work so the present choice is only a first step towards 
physical/perceptual accounting in wind instruments design.  \\
Notice also that the complex question of the nonlinear excitation 
mechanism, as the reed of an oboe or the lips of a brass player e.g., 
is discarded (equation~(\ref{eq.sturmliouville}) is homogeneous), although nonlinear effects 
can be important due to the nature of fluid dynamics, e.g. inside a reed 
or at the excitation localization in a flute. The hypothesis is made here that the nonlinear effects are, 
in a first approximation, encoded in the Fourier coefficients of the potential (see section~\ref{criterion}), 
corresponding to an arbitrary and unknown excitation. 
This is surely not satisfactory from a physical point 
of view and will have to be modified carefully in the future. 

In the present context, where losses are assumed to vanish or to be included into the boundary 
conditions, the energy inside the duct, noted $E$, is a conserved quantity~: 
\begin{eqnarray}
E=T(x,t)+U(x,t)=\frac{\pi\rho_0}{8}\int_0^LD^2(|\phi_x|^2+|\frac{\phi_t}{c}|^2) dx
\label{eq.energy}
\end{eqnarray}
with $\frac{\partial E}{\partial t}=0$ and $T(x,t), U(x,t)$ are computed 
in the appendix. 
Following the above definition of an oscillation regime, the energy can be taken as a design criterion for the duct shape. 
As it is a function of the diameter, this one can be considered as an unknown function, provided 
convenient data are given for the design. In that respect, 
as musical instruments are designed to operate in harmonic regimes, 
the potential, being in that case periodic, is amenable to a Fourier series analysis and can thus be 
written as a linear combination of time-harmonic elementary components, 
the coefficients of which will be 
other unknowns of the design problems, while the frequencies will be considered as the data. 
This will be made more precise in section~\ref{criterion}. Thanks to the computations in the appendix, 
the link is established between an oscillation regime in the sense of A.H. Benade and 
the way the Webster horn is derived, as both rely upon 
the same definite integral, namely the lagrangian action, through a Legendre 
transformation. As optimal control theory is variational too, one has a coherent set of tools for 
the design methodology to be exposed.
\subsection{\underline{The wave horn equation as a control system}}
\label{optimalcontrol}
Suppose a dynamical system is given as~:
\begin{equation}
X^{'}=f(X,U)\\
\label{eq.stateg}
\end{equation}
where $X(x)\in {\mathbb R}^n$ is the state, $U(x)\in [U_1,U_2]\subset {\mathbb R}$ the control 
and a prime still denotes the derivative 
with respect to the independent variable $x$. The control $U$ is purposedly restricted to be 
a scalar as this suffices for the needs of the present work (see~\cite{pontriaguine:1978,clarke:1983} for a more general exposition). 
The purpose of Optimal Control theory 
is to find a control $U$ such that some functional (a scalar function of the state and the control)~:
\begin{equation}
J(X,U)=\int_{x_0}^{x_1}g(X,U)dx
\label{eq.criteriong}
\end{equation}
is maximized or minimized, making it an optimization problem under dynamical constraints, named for this 
reason infinite-dimensional. This variational-type problem is solved the following way~: 
using a vector Lagrange multiplier 
$\mu$, with components $\mu_i, i=1,\ldots,n$, define the {\em Hamiltonian function} of the problem, 
a scalar function, as~: 
\begin{equation}
H(X,U,\mu)= g(X,U)+\mu^T f(X,U)
\label{eq.hamiltoniang}
\end{equation}
One then obtain {\em necessary} conditions for an optimum as the differential-algebraic system of equations (DAE)~:
\begin{equation}
\left\{\begin{array}{rcl}
X^{'}&=&\frac{\partial H}{\partial \mu}\vspace*{.1cm}\\
-\mu^{'}&=&\frac{\partial H}{\partial X}\vspace*{.1cm}\\
0&=&\frac{\partial H}{\partial U}
\end{array}\right.
\label{eq.necessaryg}
\end{equation}
together with suitable initial/boundary conditions. Notice first that, as $H$ is a scalar, all of the three 
above equations are {\em vector} ones~: $\frac{\partial H}{\partial X}$ e.g. is the vector with 
components $\frac{\partial H}{\partial X_i}, i=1,\ldots,n$. One has then a set of $2n+1$ equations, the first $2n$ 
ones being differential equations, the last one being of order zero. 
Notice also that these conditions~(\ref{eq.necessaryg}) 
are only {\em necessary} so that 
more investigations are needed to get sufficient conditions, i.e. a complete and unique solution. 
Also, the so-called Pontryagin 
maximum principle~\cite{pontriaguine:1978,clarke:1983} allows to have more precise results for the 
above when the control has e.g. bound 
constraints, as will be the case for the acoustic design problem here.\\
Now, in order to fit the acoustic design problem at hand to this control theoretic setting, 
it is useful to reformulate the second order 
equation~(\ref{eq.sturmliouville}) as a two dimensional first order 
system, suitable for control purposes in view. To this end, define the following 
variables~: $X_1=D, Y=\phi, Z=\phi^{'}(=v)$ and $W$ the vector with components $X_1,Y,Z$. Also, as the section diameter 
$D$ of the horn is an unknown to be determined as a function of $x$, it can be considered as a {control variable} 
to be designed. Actually, one can control either $D$ itself or 
the way it varies along the $x$ axis, that is one can control its derivative $D^{'}$ with respect to $x$~: as 
the modelling hypothesis i) above is that the section 
is a slowly varying function of $x$, one natural choice for the control is the derivative $D^{'}$ and bounds will have to be imposed to it in order 
to satisfy the hypothesis given by inequality~(\ref{eq.inequal}). 
As a consequence, defining the control variable $U$ as $D^{'}$, constraints on it will be {bound constraints}~: 
$D_1 \leq D^{'} \leq D_2$. 
Together with those on the control, 
constraints can also be imposed to the state through design/building constraints on the diameter $D$ itself~: 
an obvious mandatory state constraint is e.g. $X_1=D>0$. 
After having developped the second order derivatives, and with these notations, equation (\ref{eq.sturmliouville}) 
rewrites as the first order differential system~:
\begin{equation}
\left\{\begin{array}{rcl}
X^{'}_1&=&U\\
Y^{'}&=&Z\\
Z^{'}&=&-2\frac{U}{X_1}Z-k^2Y
\end{array}\right.
\label{eq.state}
\end{equation}
which is a dynamical control system, affine in the control $U$~: $W^{'}=f(W,U)=h_1(W)+h_2(W) U$, 
with immediate definitions for $h_1$ and $h_2$. One has even a {linear drift}~: $h_1(W)=A W$, where 
matrix $A$ is clearly defined. 
The singularity in $X_1=0$ is not a real problem as it corresponds to a vanishing diameter, 
a highly uninteresting situation (except possibly at one boundary, the apex of a complete cone e.g.). Now, referring to theorem I, item 2 of the appendix, for each 
eigenvalue $\lambda_n=k_n^2$ the eigenfunction $\varphi_n$ satisfies the wave 
horn equation with the same boundary conditions. One can then rewrite for each eigenvalue $\lambda_n$, the 
corresponding first order differential system of equations~:
\begin{equation}
\left\{\begin{array}{rcl}
X^{'}_{2n}&=&X_{2n+1}\\
X^{'}_{2n+1}&=&-2\frac{U}{X_1}X_{2n+1}-k_n^2X_{2n}
\end{array}\right.
\label{eq.statej}
\end{equation}
with~: $X_{2n}=\varphi_n,X_{2n+1}=\varphi_n^{'}$. This set of equations, $i=1,\ldots, $ together with the equation 
$X^{'}_1=U$ constitutes the controlled dynamical model that will enter the optimization process as constraints. 
Remark that on a practical side, only a finite number of eigenvalues will be retained, as is detailed below in 
section~\ref{criterion}.
One should notice also that, whereas in most control theoretic situations 
one looks for a {feedback control}, the searched after control for the above horn equation is intrinsically {open-loop}, as the duct is designed once for all (at least in 
the present state of technology...) so that optimal control theory and the Pontryagin Maximum Principle~\cite{pontriaguine:1978,clarke:1983} are well suited in the present context, although giving only necessary conditions, as already mentionned.
\section{Oscillation regime as an optimization criterion}
\label{criterion}
Besides the above simple, albeit new, problem rewriting, the crux in the followed approach 
of wind instrument design is in formulating the oscillation regime, as defined in section~\ref{hornoscillation}. 
It must be done in a way suitable for conducting the optimization while having the controlled dynamical model 
of section~\ref{optimalcontrol} above at hand. 
Considering the above definition for an oscillation regime, 
the potential $\phi$ can be assumed to be a time-periodic function. 
For harmonicity requirements of the signals inside the duct, the wave numbers $k_n$s' are fixed data given as multiples 
of a fundamental frequency~: $k_n=n k_0$, $k_0=\frac{2 \pi}{c}f_0$. To each of the $k_n$s' 
corresponds a solution $\varphi_n(x)$ of the Webster equation with $k_n$ as wave number. 
The set $\{\varphi_n(x)\}_{n\in{\mathbb Z}}$ is a basis of the Hilbert space of square-integrable functions 
on $[0,L]$, thanks to Theorem 1, item 3, appendix. Thus, thanks to Theorem 1, item 4, appendix, the overall potential 
$\phi$ can be written as~:
\begin{equation}
\phi(x,t)=\sum_{n\in{\mathbb Z}}c_n\varphi_ne^{i\omega_n t}
\label{eq.potential}
\end{equation}
which is seen to be the Fourier series decomposition of $\phi(x,t)$ with respect to the time variable. 
Thus Parseval theorem allows to write~:
\begin{equation}
|\phi|^2=\sum_{n\in{\mathbb Z}}c_n^2|\varphi_n|^2,|\phi^{'}|^2=\sum_{n\in{\mathbb Z}}c_n^2|\varphi_n^{'}|^2
\label{eq.squares}
\end{equation}
The energy $E$ can now be computed as a function of $c_n,\varphi_n,\varphi_n^{'}$, simply substituting for $\phi$ from 
(\ref{eq.potential}) into (\ref{eq.energy}), 
which gives, thanks to the relations~(\ref{eq.squares})~:
\begin{equation}
E=\frac{\pi\rho_0}{8}\int_0^L\left(D^2\sum_{n\in{\mathbb Z}}c_n^2(|\varphi_n^{'}|^2+\frac{\omega_n^2}{c^2}|\varphi_n|^2)\right)dx = \frac{\pi\rho_0}{8}\int_0^L\left(D^2\sum_{n\in{\mathbb Z}}c_n^2(|X_{2n+1}|^2+k_n^2|X_{2n}|^2)\right)dx 
\label{eq.energy.1}
\end{equation}
the last equality stemming from the definition of $X_{2n},X_{2n+1}$ 
at the end of section~\ref{preliminaries}.
Now assume that the energy is conserved inside the horn. At least, one can consider in a first approximation that the small part of energy which is radiated outside the duct 
or dissipated at the inside boundary is exactly compensated with the energy brought in by the excitation mechanism such as the reed of the instrument. Such dissipation and other phenomena should surely be considered in future work, possibly 
through the boundary conditions. 
The definition of an oscillation regime leads then to maximize $E$ as 
defined in equation~(\ref{eq.energy.1}). 
Thus, with {\em the $\omega_n$s'  known and fixed as given data}, the optimization problem above has 
to {\em find $D(x)$ and the $c_n$s'}. Observe nevertheless that obtaining optimal $c_n$s' does not give 
any indication on some perceptive quality factor which is indeed important for a high quality design. 
Thus it is likely at this stage that optimizing the above-defined oscillation regime is unsufficient 
to that end. 
Section~\ref{numeric} will present simulations for the above defined problem, after the design problem itself 
has been mathematically posed in the following section. 
\section{Optimal design of a horn shape}
\label{webster}
Thanks to the results of section~\ref{preliminaries} and~\ref{criterion}, the design of 
axisymmetric wind instruments can be now formulated as an optimal control problem in 
the following way. Consider an 
oscillation regime as defined and developped in sections~\ref{hornoscillation}, \ref{criterion}, 
each component being governed by the Webster 
horn equation with its own wave number $k_n$ and keeping the approximation of the 
solution to a finite number, $N$, of its first terms. 
%
%
The data are made of the set of fixed multiples of a fundamental frequency $f_0$, leading to the set 
$k_n=j_nk_0, j_n=1,\ldots, N$, with $j_n$ allowing to take into account that a complete series of 
harmonics of $f_0$ or only the odd part of this series can be present.
As $\phi$ is unknown, the $c_n$s' of its time-harmonic decomposition are unknown parameters 
and considered as constant controls. In the following, the $c_n$s' 
will be gathered in a parameter vector noted $C=(c_1,c_2,\ldots, c_N)$.
%
%
Then the searched after shape of the duct and the $c_n$s' are obtained as the solution of the following problem~:
\begin{flushleft}(${\cal P}$) {\bf Maximize} $E$ as given by equation~(\ref{eq.energy.1}), with respect to 
$D$ and $C$ \end{flushleft}
\hspace*{.6cm}{\bf subject to}~:
\begin{enumerate}
\item The $2N+1$-dimensional dynamical model $(n=1,\ldots ,N)$~:
\begin{equation}
\left\{\begin{array}{lcl}
X^{'}_1&=&U\\
X^{'}_{2n}&=&X_{2n+1} \\
X^{'}_{2n+1}&=&-2\frac{U}{X_1}X_{2n+1}-k_n^2X_{2n}
\end{array}\right.
\label{eq.states}
\end{equation}
each $n^{th}$ pair of the last $2N$ equations being a Webster equation with wave number $k_n$. 
\item Bound constraints on the variable control~: $D_1\leq U=D^{'} \leq D_2$.
\item State constraint~: $X_1\geq a> 0$, $a$ a fixed real number.
\end{enumerate}
It is apparent that, whenever $U>0, X_1(0)>0$, the state constraint is satisfied. 
Within this context, the strong Pontryagin 
Maximum Principle (PMP)~\cite{pontriaguine:1978,clarke:1983} is applicable. 
Necessary conditions for an optimal control (i.e. an 
optimal duct shape) are obtained in the following way. 
Along the general setting given by~(\ref{eq.stateg}), (\ref{eq.criteriong}), (\ref{eq.hamiltoniang}), (\ref{eq.necessaryg}), adjoin the dynamical constraints~(\ref{eq.states}) to the criterion $E$, (\ref{eq.energy.1}), through Lagrange multipliers 
$\mu=(\mu_n)_{n=1,\ldots, N}$ 
, define the Hamiltonian of problem (${\cal P}$) as~:
\begin{equation}\begin{array}{c}
H(X,U,C,\mu)=\frac{\pi\rho_0}{8}X_1^2\sum_{n=1}^Nc_n^2(X_{2n+1}^2+k_n^2X_{2n}^2)+\mu_1 U \vspace*{.5cm}+ \sum_{n=1}^N(\mu_{2n}X_{2n+1} - \mu_{2n+1}(2\frac{U}{X_1}X_{2n+1}+k_n^2X_{2n}))
\end{array}
\label{eq.hamiltonian}
\end{equation}
Then, following~(\ref{eq.necessaryg}), $\mu$ is the solution of the adjoint differential system $(n=1,\ldots ,N)$~:
\begin{equation}
\left\{\begin{array}{lcl}
\mu_{1}^{'}&=& -\frac{\partial H}{\partial X_1}=-\frac{\pi\rho_0}{4}X_1\sum_{n=1}^Nc_n^2(X_{2n+1}^2+k_n^2X_{2n}^2)-2\frac{U}{X_1^2}\sum_{n=1}^N\mu_{2n+1}X_{2n+1}\vspace*{.2cm}\\
\mu_{2n}^{'}&=& -\frac{\partial H}{\partial X_{2n}}=k_n^2(\mu_{2n+1}-\frac{\pi\rho_0}{4}c_n^2X_1^2X_{2n})\vspace*{.2cm}\\
\mu_{2n+1}^{'}&=& -\frac{\partial H}{\partial X_{2n+1}}=2\frac{U}{X_1}\mu_{2n+1}-\mu_{2n}-\frac{\pi\rho_0}{4}c_n^2X_1^2X_{2n+1}
\end{array}\right.
\label{eq.costate}
\end{equation}
Thus the optimization problem will be complete when initial-boundary conditions are 
specified. Notice that time-initial conditions are not relevant here as the time variable has been eliminated 
through the hypothesis of time-harmonic regime and of energy conservation~: the outputs of the 
optimization process are space-dependent, giving e.g. the mode shapes inside 
the duct. As a starting point and in order to simply illustrate the method, without going into deep physical 
considerations of a specific instrument, 
that are postponed to a forthcoming paper, the 
end at $x=0$ is assumed to be closed (typically where the excitation would be placed) so that the velocity and 
thus $\varphi^{'}$ 
vanish there~: $\forall n=1,\ldots, N, \varphi_{2n}^{'}(0)=X_{2n+1}(0)=0$. Also the other end is 
assumed to be open (e.g. the other end of a brass instrument) so that, at the first order approximation, 
the pressure and thus $\varphi$ vanish~: 
$\forall n=1,\ldots, N, \varphi_{2n}(L)=X_{2n}(L)=0$. These conditions lead to a Two-Point Boundary 
Value Problem (TPBVP) for the differential system~(\ref{eq.states}),(\ref{eq.costate})~: 
one part of the state is fixed at one end 
and another part is fixed at the other end. 
When some state component $X_j$ is left unspecified at one end, the corresponding costate component $\mu_j$ must 
vanish there~\cite{pontriaguine:1978}. Thus in the present situation~: $\mu_{2n}(0)=\mu_{2n+1}(L)=0$.
Other conditions can be imposed with more realistic considerations. 
For computing the optimal control, the Pontryagin maximum principle 
implies that, for $\hat{X}$ (resp. $\hat{\mu}$) 
a solution of the state (resp. adjoint state) equation, an optimal control 
$\hat{U}$ and optimal parameter vector $\hat{C}$ are such that~:
\begin{equation}
\hat{H}=H(\hat{X},\hat{U},\hat{C},\hat{\mu})=\max_{\{D_1\leq U\leq D_2],C=Cst\}}H(\hat{X},U,C,\hat{\mu})
\label{eq.hamiltonianopt}
\end{equation}
But one can observe that $H$ is affine with respect to $U$ so that~:
\begin{equation}
\frac{\partial H}{\partial U}=\mu_1 -\frac{2}{X_1}\sum_{n=1}^N \mu_{2n+1}X_{2n+1}=0
\label{eq.Hu}
\end{equation}
which is also named the {\em switching function} for this problem as its zeroes and its derivative 
can help determine the shape of the duct. 
As equation~(\ref{eq.Hu}) does not allow to compute the control $U$ explicitely, one is faced with 
a problem with {\em singular extremal arcs}~\cite{bryson:75}, chap.~8.
The standard method is to compute successive derivatives of $\partial H/\partial U$ 
with respect to $x$ until one is able 
to get an expression for $U$. In the present case, two such derivatives allow this, after 
tedious but straightforward computations. In case $U$ is bounded, as it is the case here, the minimum 
of the hamiltonian $H$ with respect to $U$, where $\partial H/\partial U=0$, can occur at the boundary of the domain 
(see~\cite{bryson:75,pontriaguine:1978}). This fact can be observed in the simulations~: 
the control can be on its bound along several 
intervals of the integration interval (figure~\ref{figure-2-3}) or can go from one bound 
to the other (figure~\ref{figure-10-3}), revealing a {\em bang-bang} type control~\cite{bryson:75}. 
The conclusion is that one generally obtains quasi-cones. Cones could not be obtained for 
exactly harmonic frequencies~\cite{kergomard:2008}, chap. 7. This is compatible with 
the results obtained in~\cite{dalmont-kergomard:1994}.
Now, thanks to the developments of this section, one has sufficient material to illustrate the design method 
through a few simulations, on which some qualitative observations can 
be made and confirmation of the theoretical results is given. Remember nevertheless that in the following, 
only candidates for optimal shapes are obtained at this stage (necessary conditions) 
because, due to theorem 2, only $q=2D^{'}/D$ can be uniquely 
determined from the given boundary spectral data. Further data and theoretical investigations 
are necessary towards a really optimal shape. 
\section{Numerical results}
\label{numeric}
\addtocounter{paragraph}{-3}
The numerical results presented here must be considered 
as qualitative illustrations of the above theoretical results. 
No interpretation in terms of musical quality will be attempted at this stage. 
Such results will be pursued elsewhere more deeply in order to derive conclusions on realistic physical and 
perceptual basis. Nevertheless, some features are worth noticing at this stage, 
when faced with some characteristics of real instruments. 
\subsection{Data}
\label{data}
In the following, data are as follows~: the air mass density is taken to be $\rho_0=1$, 
the sound velocity is that in free space at standard temperature~: $c=340 m/s$ and one 
focuses on a fundamental frequency for the note usually labelled $A_4$, i.e. $f=440 Hz$, leading to the 
wavelength $\lambda=c/f=0.772 m$. The chosen wave numbers are then indicated under each figure~: 
$k_1=2\pi/\lambda, k_i=n_i k_1, i=1,\ldots,n$ and $n_i$ an even or odd integer. 
At the narrow end, $x=0$, one has $D(0)$ fixed with different values for each 
simulation (see below) and $\phi^{'}(0)=0$. 
At $x=L$, the condition $\phi(0)=0$ 
is satisfied only approximately, by adjunction to the criterion through a simple penalty method. 
The duct diameter derivative is allowed to vary in the bounded 
interval $[D^{'}_m, D^{'}_M]$, $D^{'}_m, D^{'}_M$ being indicated on each figure. They are chosen to satisfy inequality~(\ref{eq.inequal}) with 
$R=D/2$, for the maximum value $k$ of the $k_n$s', 
i.e. are within the validity domain of the wave horn equation, according to~(\ref{eq.inequal}).
The interval of integration, i.e. the length of the duct, is taken to be $L=\lambda$ but could be taken as 
an unknown as well, leading to an analogous of the so-called {\em time-optimal} control. 
Last, all unspecified initial values for the state variables and all initial controls are taken to be 
random variables uniformly distributed on $[0,1]$. 
\subsection{Numerical method}
\label{method}
The numerical method used to solve the optimization problem exposed in section~\ref{webster} is the following 
(see e.g.~\cite{bryson:75}, chap.~7 for more details on several possible numerical methods used 
for optimal control problems)~: 
an initial guess is given for the control vector and, using this control, 
the state equation is numerically integrated from $X(0)$, considering too the unspecified initial conditions 
as constant controls, which are to be determined with the optimization process. 
The integration was performed thanks to a predictor-corrector scheme with an explicit Euler 
method for the prediction step and a Crank-Nicholson scheme for the correction. 
Then the costate equation is integrated backwards -notice that this equation is always linear in the costate- 
using the obtained terminal values for $X(L)$. 
This allows to compute the gradient of the objective function with respect to the control, thanks to 
$\partial H/\partial U$~\cite{clarke:1983}. 
Then a standard optimization routine is used, passing through the previous steps 
in a recursive way, in order to make this gradient decrease, 
until the specified terminal conditions at $x=L$ are -approximately here- achieved and the criterion evolves x
no more, to a specified precision ($10^{-5}$). For the time being, the criterion and the terminal constraints have 
been gathered in a single objective function through a simple penalty method but this could be improved. 
For the optimization, a quasi-Newton 
method with projection has been used together with a BFGS method for an estimation of the hessian. 
Such an approach is known to have, as usually Newton-type methods, the drawback that it can give local 
minima and to be sensitive to the initial guess for the solution~\cite{dennis:1983}. Thus, other methods such as direct 
solution of the TPBVP by shooting or multiple shooting techniques should be interesting to investigate but this 
is deferred to future work because an extra work has to be done in order to have an explicit 
expression for the control, as a function of the state and costate (see the discussion in section~\ref{webster} above). 
Nevertheless several different initial conditions for the state 
$X$ were tested : the results, not exposed here, did not show to be that sensitive but this should be confirmed 
theoretically. 
\subsection{Simulations}
\label{results}
The few simulations presented here have been done while fixing, respectively, two, 
five and ten components for the overall potential 
inside the duct, i.e. $N=2, 5, 10$ in equation~(\ref{eq.states}). In each case, the duct shape is shown 
first, followed by the corresponding modes shapes that have been normalized to their maximum value, at $x=0$. 
In addition for the case of two components only, the diameter derivative has been shown to illustrate 
the fact that it can be only piecewise continuous and that the phenomenon, mentionned in section~\ref{webster}, 
of singular extremal arcs joining regular arcs can appear (see figure~\ref{figure-2-3})~: the control 
$U=D^{'}$ is on its bound on some subintervals. 
For the five components cases, the shape is also quasi-conic but in a less obvious manner.
But for the ten component case, one has an example of abrupt change in conicity and $D^{'}$ in 
figure~\ref{figure-10-3} illustrates the possibility of bang-bang type control mentionned in section~\ref{webster}. 
On each figure the imposed lower ($D_1$) and upper ($D_2$) bounds are indicated. 
Notice that the impedance, $p/u$ and the instantaneous power, $pu$ ($p$, the pressure, $u$ the volume velocity), 
for all $x$ are easily obtained from 
the outputs, $\varphi(x), \varphi^{'}(x)$ and $D(x)$ $\forall x\in [0, L]$, of the optimization process. 
The simulations are as follows~: 
\begin{enumerate}
\item \underline{Two components}~: 
\label{FixTwo}
For this first simulation, the fundamental frequency is $A_4$ and the second 
component is the second harmonic (double frequency). The diameter at the 
closed end is $D(0)=2cm$. 
The duct shape is shown in figure~\ref{figure-2-1}, the two modal shapes in figure~\ref{figure-2-2}. 
$D^{'}$ is shown in figure~\ref{figure-2-3} as it makes appear the phenomenon described in section~\ref{webster}, 
where the control is at the bound on two subintervals, with the consequence that one obtains a quasi-cone for the duct, 
made of conic pieces joined by more complex but smooth parts.
\item \underline{Five components}~: 
\label{FixFive}
In this second simulation (see figures~\ref{figure-5-1},\ref{figure-5-2}), 
the data are the same as for two components, except for $D(0)=1cm$. 
The five components have 
frequencies $f_i=(i+1)*f_0 ; i=0,\ldots,4$. One has not shown $D^{'}$ here but the result is here again a quasi-cone, 
although this is not apparent again on figure~\ref{figure-5-1}, because $D^{'}$ varies much less but still varies 
along the duct.
\item \underline{Ten components}~: 
\label{FixTen}
In this third simulation (see figures~\ref{figure-10-1},\ref{figure-10-2},\ref{figure-10-3}), 
the data is the same as for five components, except for $D(0)=5cm$. 
The ten components have 
frequencies $f_i=(i+1)*f_0 ; i=0,\ldots,9$.  The results are noticeable as 
here again, three conic pieces are found so that the derivative $D^{'}$ is shown 
in figure~\ref{figure-10-3}, making it here again a quasi-cone, with a bang-bang type behaviour for 
$D^{'}$. The duct is made of conic pieces with different conicity joined together. Only the five first modal shapes 
have been displayed, for better readability. 
\end{enumerate}
\subsection{Qualitative observations and commentaries}
The above simulations show only qualitative albeit important results at the present stage. The main one is that, 
as foreseen by the theoretical investigations, and for the specific chosen criterion, the duct shapes are 
piecewise continuous, this being particularly clear on the case of two and ten components ; 
observe too in both cases the fact 
that the control is at the bound on some subintervals, meaning that regular arcs and singular arcs can coexist. 
One can see also that quasi-cones are obtained, in an obvious way for the first and second simulations. 
Recall nevertheless that, as mentionned, the numerical 
method can be sensitive to initial conditions in the optimization 
process, so that it is difficult to interpret these results in a precise 
manner at the present stage. The qualitative observation is that the shapes 
are quasi-cones. At the same time, the modal shapes (i.e. the 
eigenfunctions of the Sturm-Liouville problem) show a behaviour that is qualitatively 
in agreement with what is theoretically predicted for cones 
(see~\cite{kergomard:2008}, chap. 7), i.e. one has modal shapes that decrease as $x$ increases to $L$.
On one hand, these results are coherent with what is known for such geometries~\cite{dalmont-kergomard:1994}. 
On another hand, the 
fact that one gets quasi-cones essentially would lend to conclude that the used notion of oscillation 
regime does not allow to capture all the important features that are sought after for musical quality, as 
flaring horns in brasses e.g. cannot be achieved with it~:  
significant and important nonlinear phenomena~\cite{kergomard:2008} 
are not taken into account with the here retained model and design criterion. Thus, using more realistic models 
and refined criteria should be used instead.\\
The computed modes have an extremum at the closed end, coherent 
with the chosen boundary condition ($\phi^{'}(0)=0$). 
At the other open end, the modes  does not vanish exactly~: one reason for this 
is linked to the way this end condition has been taken into account in the numerical method~: a simple penalty 
method makes a compromise between satisfying the criterion and the boundary conditions there. The 
above suggested method of directly solving the TPBVP through shooting techniques should give better results 
from this point of view.
\section{Conclusion and perspectives}
\label{conclusion}
Thanks to a reformulation within a control theory framework, the design of axisymmetric wind 
instruments has been revisited. It has been shown that in order for an oscillation regime in the sense 
of A.H. Benade to take place inside the air column and with the considered linear model, the shape 
tends to be conic, as the number of time-harmonic components grows. This tends to show that such 
a criterion, as it has been mathematically described, is unsufficient to grasp the necessary 
conditions that lead to high-quality instruments. No other conclusion 
from this point of view is given at this stage.
Nevertheless, and this was one main purpose of this work, 
the approach makes it very flexible to deal with a great variety of design constraints, 
either on the control or on the shape itself. 
At the present time, only qualitative results have been given to illustrate the theoretical results 
and the field is wide open to numerical investigations as well as 
to experiments on more physical and musical premisses. In that respect, some important issues can 
be straightforwardly put in perspective~: 1) A first possible development is to use the method in a model 
fitting way, using an experimentally 
measured input impedance from a real instrument and taking a distance measure with the corresponding model 
output as criterion, in a similar way as in~\cite{kausel:2001,noreland:2003-2}. More realistic 
duct shapes are to be expected and it will be interesting to compare results from different approaches. 
2)~Refine the question of physical/perceptual 
criteria and of the horn model~: this relies upon physical considerations as well as on experiments 
with real existing instruments played by expert musicians. This is surely 
a central point to get a practical tool useful for real design. One consequence of modifying the 
criterion is that the structure of the obtained shape will likely be different from that obtained 
in the present study.
3) Develop the method for woodwinds, thereby including toneholes in the design process. 
4) Take into account some imposed-shaped parts such as pieces of cylinders 
or cones from built-in components, such as for brasses. This implies using {\em path constraints} in the 
optimal control framework. 5) All the previous points depend on the modelling question of 
how to include the nonlinear excitation 
mechanism (reed, lips...) in the model as well as other nonlinear 
phenomena, which were discarded in the present study. 6) Last, an interesting question is to try to retrieve, 
with the presented approach, the lattice of sound tubes results found in~\cite{dalmont-kergomard:1994}. This will imply to 
modify the mathematical framework of the controlled dynamical system, i.e. the function 
spaces where the state and the control 
live, in order to be able to take into account discontinuous, and not only piecewise continuous, duct shapes. 
The mathematical theory of distributions could be an adequate framework for this. 
These items are but a small part of what should be addressed within the proposed framework to achieve 
a satisfying design methodology. 
\section{Acknowledgements}
The author is undebted to two anonymous reviewers an to the editor of the special issue, J. Kergormard, 
for their numerous advices that helped greatly improve the paper.

\appendix
\label{appendix}
\paragraph{\underline{The horn wave equation}}
\label{hornequation} 
Instead of the usual, Newton-type modelling approach, making use 
of forces and moments balance equations, a lagrangian one is 
followed here, as it fits better the variational context for the design 
methods that form the core of the present work.
For standard derivation the reader is referred to~\cite{putland:1994,bruneau:1998} for 
physical considerations and to~\cite{levey-isma07} for control theory formulation.\\
For a general pressure field, the lagrangian action density at each 
time instant inside the duct is the difference between a kinetic term and a 
potential term, that are computed in a standard way, as follows. Assumption 
ii) implies that the duct can be considered as a continuous 
stack of cross-sections $S(x)$, parametrized by the abcissa $x$. For each 
section $S(x)$ located at $x$ along the horn axis, 
an action density is computed as the integral of densities of the 
particles over the section. This leads to 
an expression proportional to the cross-section area, i.e. to $D^2(x)$. 
Firstly, the kinetic term writes~: 
\begin{equation}
T(x,t)=\int_{S(x)}\frac{1}{2}\rho_0v^2d\sigma=\frac{\pi D^2}{4}\frac{\rho_0}{2}v^2=\frac{\pi D^2}{4}\frac{\rho_0}{2}|\phi_x|^2
\label{eq.kinetic}
\end{equation}
Similarly, the potential energy term is given as~: 
\begin{equation}
U(x,t)=\int_{S(x)}\frac{p^2}{2\rho_0 c^2}d\sigma=\frac{\pi D^2}{4}\frac{p^2}{2\rho_0 c^2}=\frac{\pi D^2}{4}\frac{\rho_0}{2}|\frac{\phi_t}{c}|^2
\label{eq.potentiel}
\end{equation}
Thus the lagrangian action density writes~:
\begin{equation}
L(x,t)=T(x,t)-U(x,t)=\frac{\pi D^2}{4}\frac{\rho_0}{2} (|\phi_x|^2-|\frac{\phi_t}{c}|^2) 
\label{eq.lagrangiandensity}
\end{equation}
Eventually, the lagrangian action inside the duct writes~:
\begin{equation}
{\cal L}=\int_{t_0}^{t_1}\int_0^L L(x,t) dx dt=\frac{\pi\rho_0}{8} \int_{t_0}^{t_1}\int_0^LD^2(|\phi_x|^2-|\frac{\phi_t}{c}|^2) dx dt
\label{eq.lagrangianaction}
\end{equation}
According to Hamilton's stationary action principle, the dynamics inside 
the duct is obtained as the Euler equation of the above action ${\cal L}$~: 
\begin{equation}
\frac{\partial L}{\partial \phi}-\frac{\partial }{\partial x}(\frac{\partial L}{\partial \phi_x})-\frac{\partial }{\partial t}(\frac{\partial L}{\partial \phi_t})=0
\label{eq.eulerapp}
\end{equation}
But one can see that $\frac{\partial L}{\partial \phi}=0$, thus~:
\begin{equation}
-\frac{\partial }{\partial x}(2D^2\phi_x)+2(D^2/c^2)\phi_{tt}=0
\label{eq.eulerapp-1}
\end{equation}
and eventually, dividing by $2D^2$~:
\begin{eqnarray}
\frac{1}{c^2}\phi_{tt}-\phi_{xx}-2\frac{D_x}{D}\phi_x=0
\end{eqnarray}
which is recognized to be the wave horn equation or Webster horn equation. 
The interest of the above variational derivation for this equation 
lies first in the fact that the optimal control approach to the design is variational in nature too. 
Thus, due to the harmonic nature of waves inside musical instruments, the horn equation writes, 
for one fixed value $\omega$~:
\begin{equation}
\phi_{xx}+2\frac{D_x}{D}\phi_x+\frac{\omega^2}{c^2}\phi=0
\label{eq.hornapp}
\end{equation}
i.e. using primes from now on to denote the derivatives with respect to the 
spatial dimension $x$ along the axis, the only remaining independent variable~:
\begin{equation}
\begin{array}{l}
\phi^{''}+2\frac{D^{'}}{D}\phi^{'}+k^2 \phi=0
\label{eq.sturmliouvilleapp}
\end{array}
\end{equation}
this equation being adjuncted a set of suitable 
boundary conditions~: 
\begin{equation}
\begin{array}{l}
a_1\phi(0)+b_1\phi^{'}(0)=0\;,\;a_2\phi(L)+b_2\phi^{'}(L)=0
\label{eq.BC}
\end{array}
\end{equation}
which represents a resonator without active components and for which losses can be taken into account in the boundary 
conditions that will be precised for the simulations in section~\ref{numeric}. 
Thus one is faced with a homogeneous Sturm-Liouville problem. 
\paragraph{\underline{Elements of Sturm-Liouville theory}}
\label{sturmliouville}
Two theorems from spectral theory of Sturm-Liouville problems 
are recalled here for self-containedness~: one~\cite{dieudonne:1960} for the so-called ``direct problems'' 
and the second~\cite{katchalov:2001} for ``inverse problems''. 
Let $I=[O,L]$ and $q(x)=2\frac{D^{'}}{D}$. Then ~: \\
{\bf Theorem 1}~\cite{dieudonne:1960}~: For every function $q(x)$ continuous in $I$~:
\begin{enumerate}
\item The Sturm-Liouville problem has an infinite strictly increasing sequence of eigenvalues $\lambda_n\in \mathbb R$ such that $\lim_{n\rightarrow\infty} \lambda_n=+\infty$ and 
the series $\sum_n 1/\lambda_n^2$ converges.
\item For each eigenvalue $\lambda_n$, the homogeneous Sturm-Liouville problem has a real-valued solution $\varphi_n(x)$ such that $\int_a^b \varphi_n^2(x) dx = 1$, which is unique up to a multiplicative real constant. 
\item The sequence $(\varphi_n)$ is an orthonormal system in a convenient Hilbert space of functions.
\item Let $w$ be complex-valued continuous function defined in $I$, the primitive of a ruled function $w^{'}$ such that~: (i) $w^{'}$ is continuous in $I$, except possibly at a 
finite number of interior points. (ii) $w^{'}$ has a derivative $w^{''}$ continuous in every interval where $w^{'}$ is continuous. (iii) $w$ satisfies the boundary conditions in~(\ref{eq.BC}).
Then, if $c_n=<w,\varphi_n>=\int_{a}^bw(s)\varphi_n(s)ds$, one has~: $w(x)=\sum_n c_n\varphi_n$ where the series converges uniformly and absolutely in $I$.
\end{enumerate}
On another hand, the design question itself relies upon the 
following inverse boundary spectral theorem~:\\
{\bf Theorem 2}~\cite{katchalov:2001}~: Assume that 
$\{\lambda_1, \lambda_2,\ldots, \varphi_1^{'}(0),\varphi_2^{'}(0),\ldots\} $ 
are the boundary spectral data of the Dirichlet-Schr\"odinger operator, ${\cal A}_0=-\frac{d^2}{dx^2}+q$, 
corresponding to the above Sturm-Liouville problem, on an interval 
$[0, L]$. Then, these data determine $L$ and $q(x)$ uniquely.\\
It is moreover worth noting that only $q$ can be determined uniquely from these data. This means 
here that only the ratio $D^{'}/D$ is so. Thus to determine $D$ itself needs supplementary data. 
%
\begin{figure}[b]
\centering
\psfrag{xx}{$x$ (in m)}
\psfrag{DUCT}{Duct shape}
\psfrag{DU}{$D$ (in m)}
\includegraphics[width=7in,height=3in]{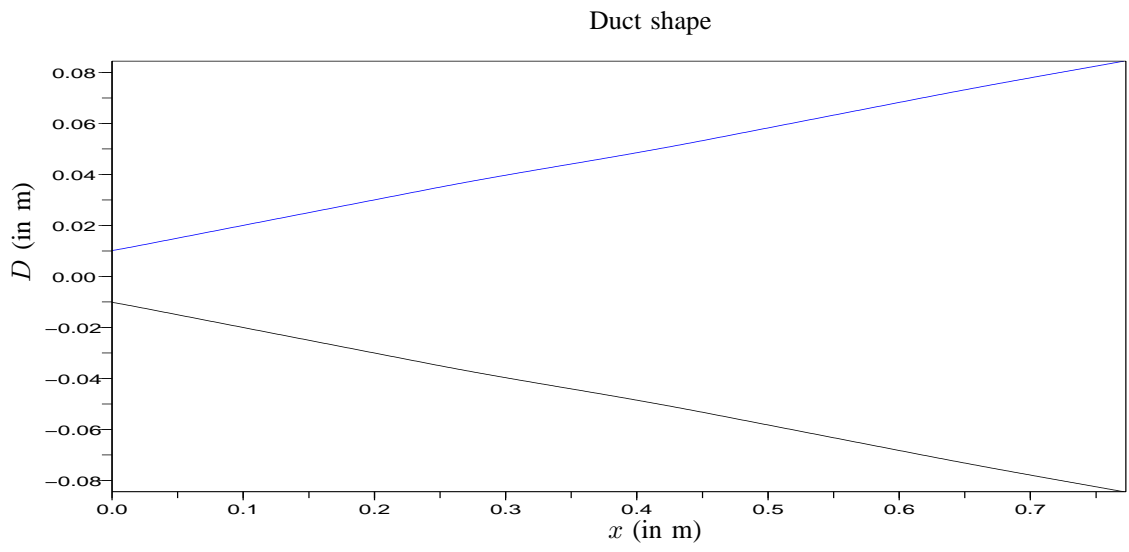}
\caption{Duct shape, 2 components, $-0.2\leq D^{'}\leq 0.2$}
\label{figure-2-1}
\end{figure}
\begin{figure}[b]
\centering
\psfrag{POT}{$\varphi_1(x),\varphi_2(x)$}
\psfrag{xx}{$x$ (in m)}
\psfrag{PHI}{$\varphi_n$}
\includegraphics[width=7in,height=3in]{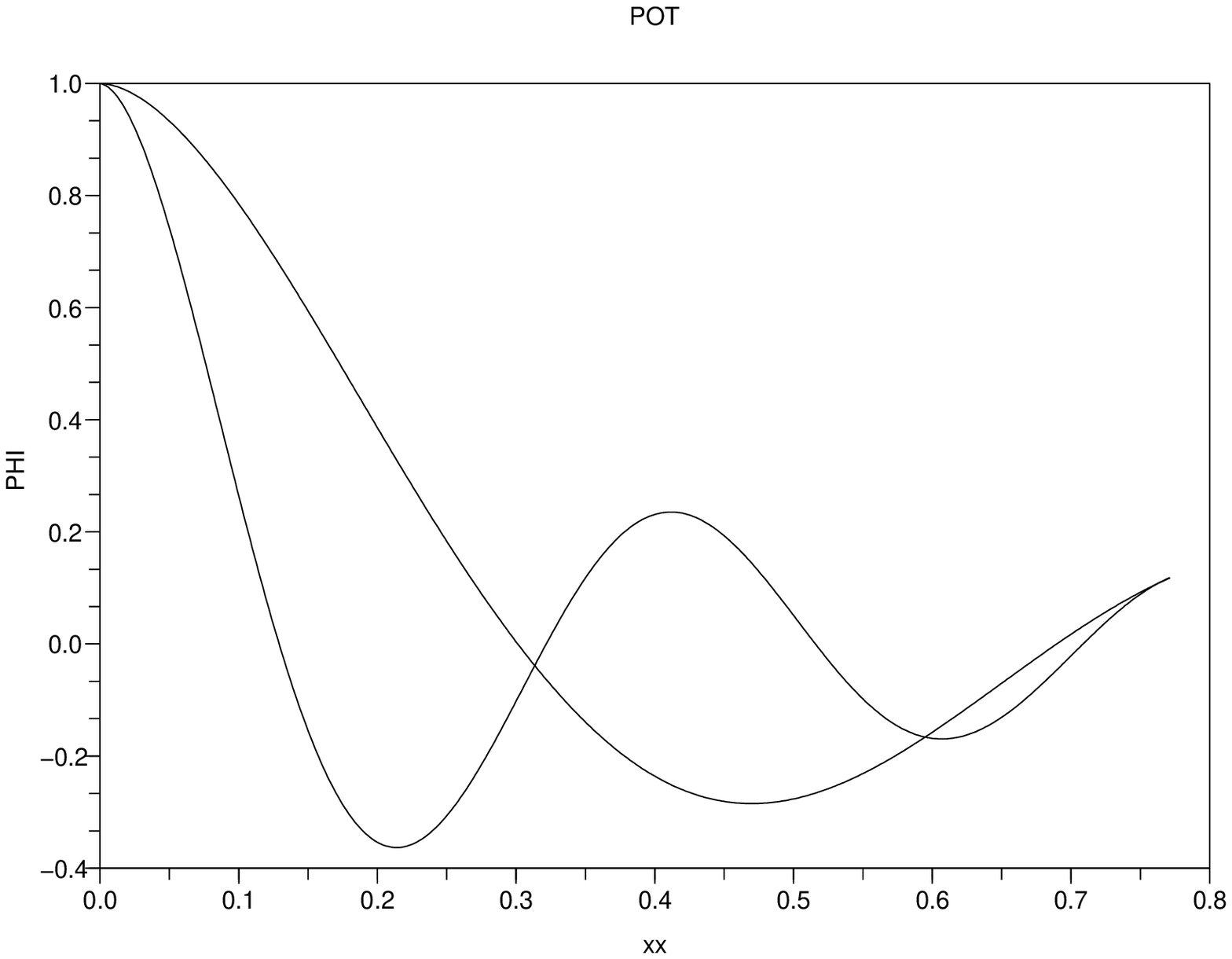}
\caption{Modal shapes, 2 components, $-0.2\leq D^{'}\leq 0.2$}
\label{figure-2-2}
\end{figure}
\begin{figure}[b]
\centering
\psfrag{xx}{$x$ (in m)}
\psfrag{DER}{$D^{'}(x)$}
\psfrag{DD}{$D^{'}$}
\includegraphics[width=7in,height=3in]{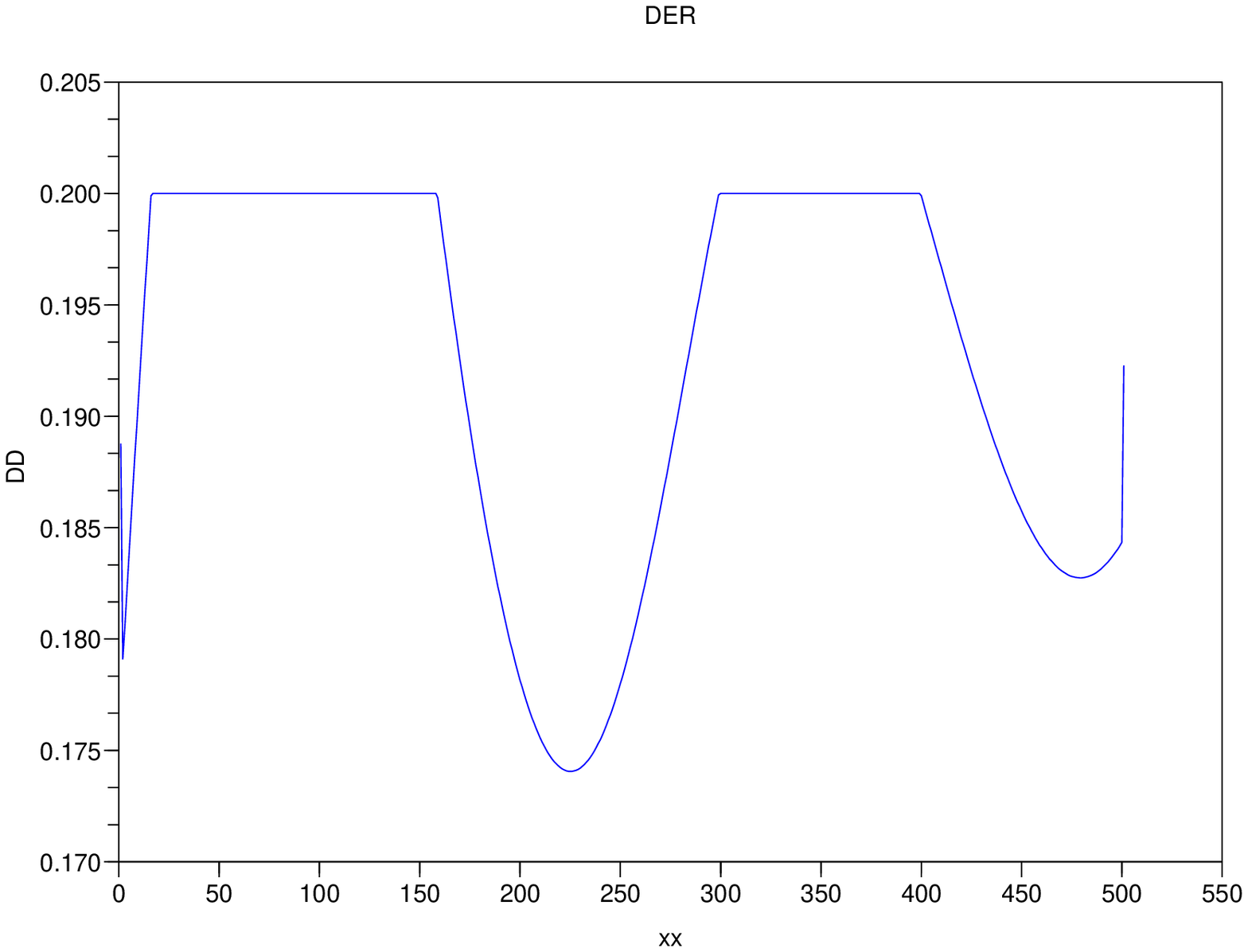}
\caption{Derivative of the duct diameter, 2 components, $-0.2\leq D^{'}\leq 0.2$}
\label{figure-2-3}
\end{figure}
%
\begin{figure}[b]
\centering
\psfrag{xx}{$x$ in m}
\psfrag{DUCT}{Duct shape}
\psfrag{DU}{$D$ in m}
\includegraphics[width=7in,height=3in]{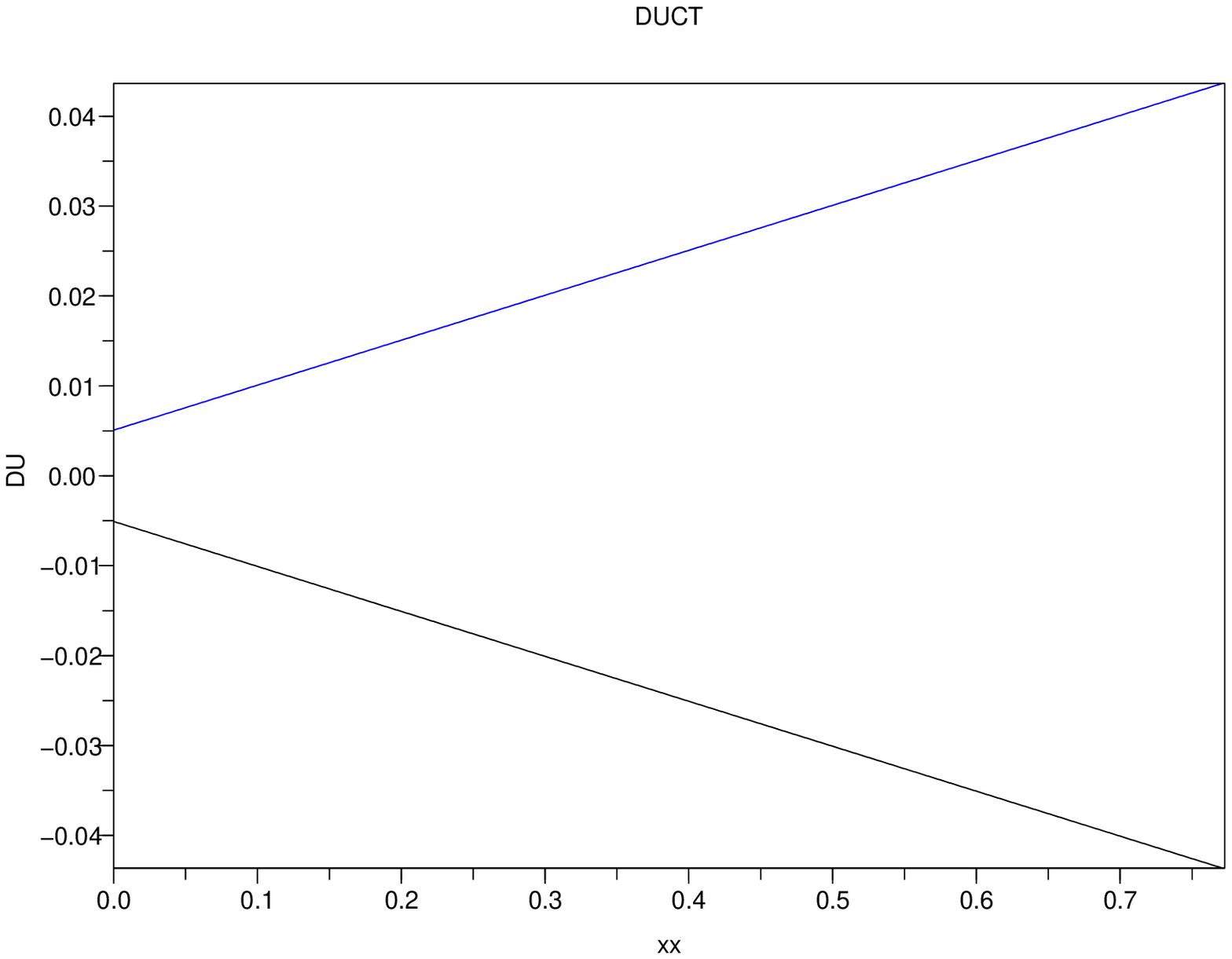}
\caption{Duct shape, 5 components, $-0.2\leq D^{'}\leq 0.2$}
\label{figure-5-1}
\end{figure}
\begin{figure}[b]
\centering
\psfrag{POT}{$\varphi_i(x), i=1,\ldots,5$}
\psfrag{xx}{$x$ (in m)}
\psfrag{PHI}{$\varphi_n$}
\includegraphics[width=7in,height=3in]{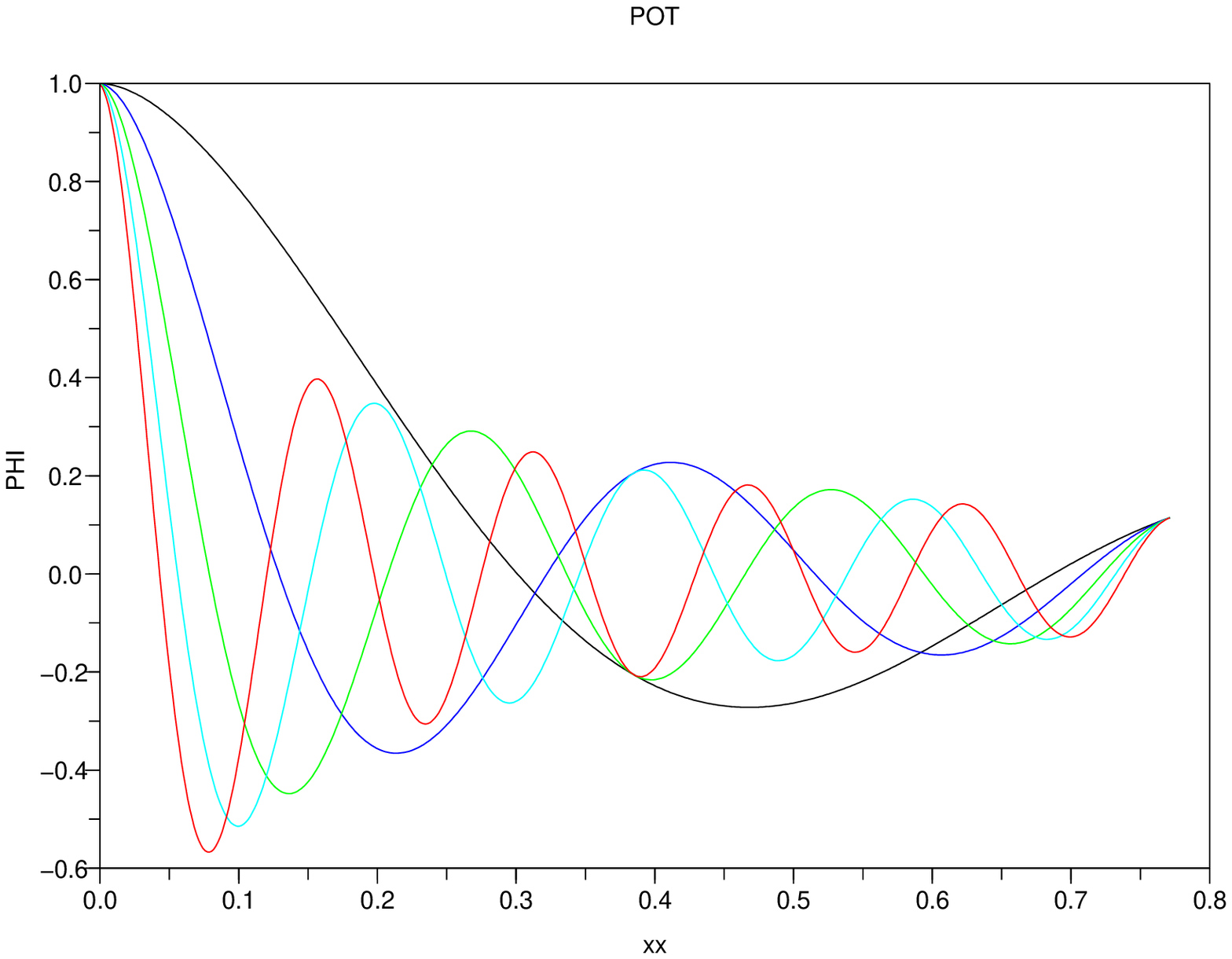}
\caption{Modal shapes, 5 components, $-0.2\leq D^{'}\leq 0.2$}
\label{figure-5-2}
\end{figure}
%
\begin{figure}[b]
\centering
\psfrag{xx}{$x$ in m}
\psfrag{DUCT}{Duct shape}
\psfrag{DU}{$D$ in m}
\includegraphics[width=7in,height=3in]{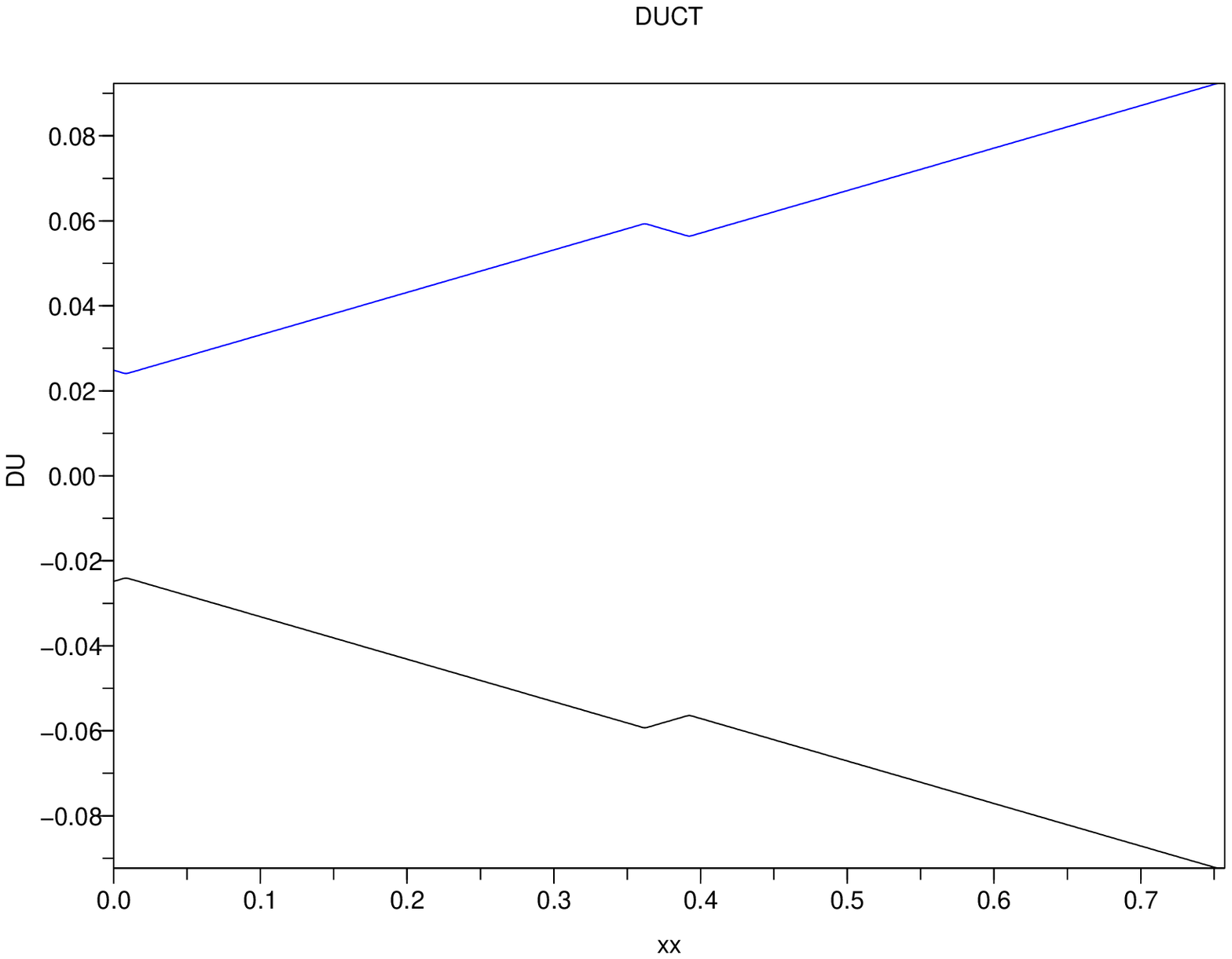}
\caption{Duct shape, 10 components, $-0.2\leq D^{'}\leq 0.2$}
\label{figure-10-1}
\end{figure}
\begin{figure}[b]
\centering
\psfrag{POT}{$\varphi_i(x), i=1,\ldots,5$}
\psfrag{xx}{$x$ (in m)}
\psfrag{PHI}{$\varphi_n$}
\includegraphics[width=7in,height=3in]{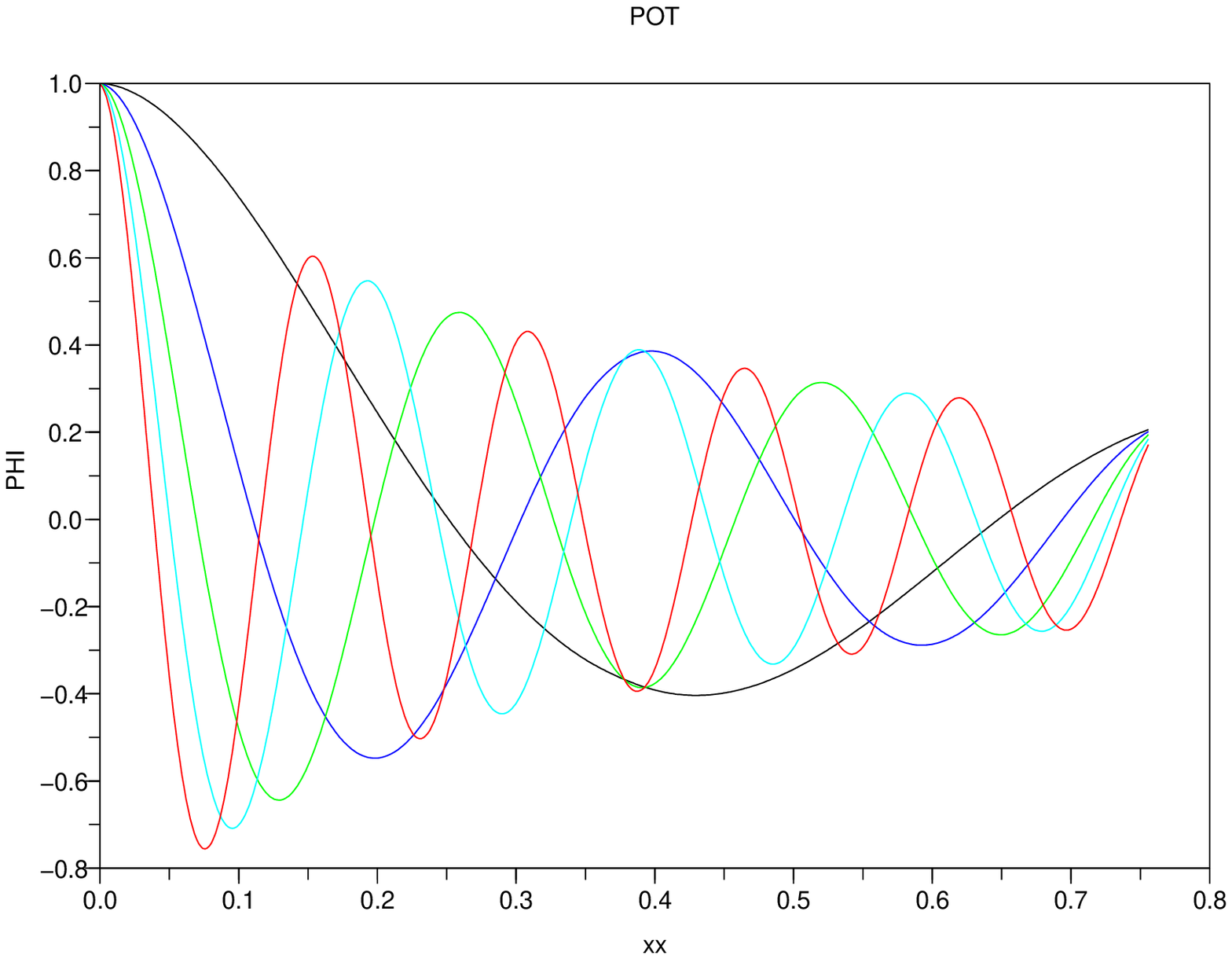}
\caption{Modal shapes, 10 components, $-0.2\leq D^{'}\leq 0.2$}
\label{figure-10-2}
\end{figure}
\begin{figure}[b]
\centering
\psfrag{xx}{$x$ (in m)}
\psfrag{DER}{$D^{'}(x)$}
\psfrag{DD}{$D^{'}$}
\includegraphics[width=7in,height=3in]{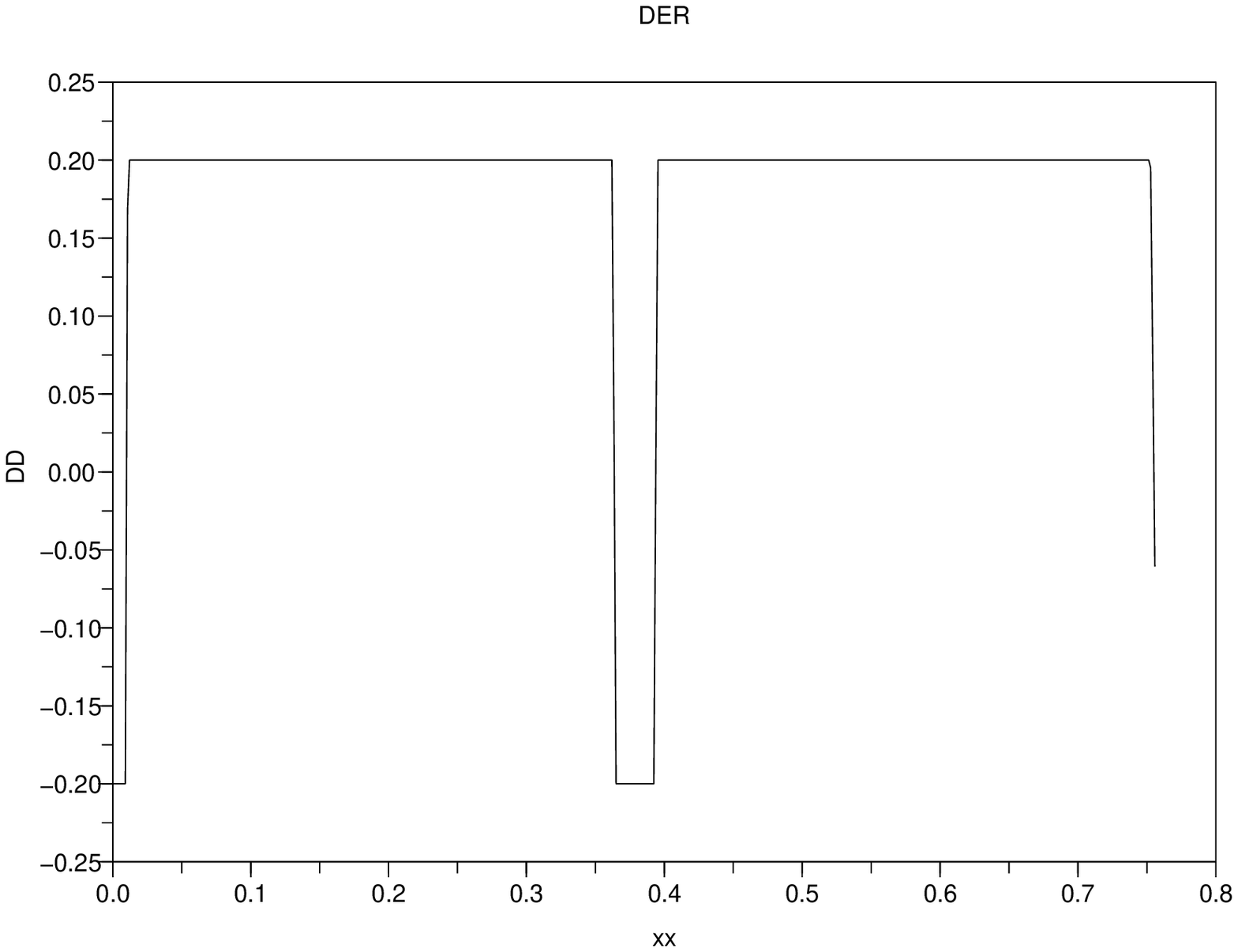}
\caption{Derivative of the duct diameter, 10 components, $-0.2\leq D^{'}\leq 0.2$}
\label{figure-10-3}
\end{figure}
\end{document}